\documentclass{amsart}
\usepackage{amsmath,amsthm,amssymb,IMjournal}

\usepackage{graphicx} 
\usepackage{amsmath}

\usepackage{mathtools}
\usepackage
{amstext, amscd, setspace, tikz-cd,mathtools, enumerate, graphics, latexsym, siunitx}

\usepackage{pst-node}
\usepackage{tikz-cd} 

\usepackage{amsthm}

\newtheorem*{mainthm}{Main Theorem}

\newcommand{\MainTheorem}{\hyperref[refer]{Main Theorem}}

\usepackage{amssymb}

\usepackage{hyperref}

\usepackage{PTSerif} 

\usepackage[cmtip,all]{xy}
\newcommand{\properideal}{%
\mathrel{\ooalign{$\lneq$\cr\raise.22ex\hbox{$\lhd$}\cr}}}

\newcommand{\properring}{%
\mathrel{\ooalign{$\gneq$\cr\raise.22ex\hbox{$\rhd$}\cr}}}

\newcommand{\GL}{{\rm GL}}

\numberwithin{equation}{section}

\catcode`\@=11
\catcode`\@=12

\theoremstyle{definition}
\numberwithin{equation}{subsubsection}

\begin{document}


\title{On Minimal generating sets of splitting field, Cluster Towers and Multiple Transitivity of Galois Groups}

\author{\|Shubham  |Jaiswal |\ORC{0009-0005-6769-2907}, 
        \|P |Vanchinathan|}



\abstract 
   A natural generating set for a Galois extension regarded as
the splitting field of  an irreducible polynomial is introduced and investigated here.
Minimal generating sets arising in this context
 throw  many surprises compared to the analogous concept in the context of vector spaces:  they can be of different cardinalities. In fact we establish that for a certain family of polynomials over the rationals, we have minimal generating sets of all cardinalities in a certain range and that these are the only possible cardinalities for minimal generating set for such a polynomial. 
 We also study how minimal generating sets behave under multiple transitivity of the Galois group and consequently prove the existence of polynomials with all minimal generating sets of uniformly same cardinality. We also connect minimal generating sets with the concept of root cluster tower of an irreducible polynomial introduced in M Krithika, P Vanchinathan (2024).
\endabstract

\keywords
minimal generating sets of splitting field of polynomial; cluster towers; multiple transitivity of Galois groups
\endkeywords

\subjclass
11R32, 12F05, 12F10, 20B20, 20B35
\endsubjclass

\section{Introduction}

Suppose for an irreducible polynomial $f$ of degree $n$ over base field $K$, the Galois group is $\mathfrak S_n$. It means every permutation of its roots gives rise to an automorphism of the splitting field. In other words knowing an automorphism at a subset of roots does not determine the automorphism fully. In this case we can say loosely that there is no interrelationship among the roots of $f$. When the Galois group is not the whole of $\mathfrak S_n$, it implicitly says there is a relationship among the roots. For example if the field generated by a certain subset of roots of $f$ contains some other roots then it is clear that Galois group is a proper subgroup of $\mathfrak S_n$. A famous result of Galois says that for a prime degree polynomial, the
splitting field is generated by two roots when the Galois group is solvable (Refer to $\S$ 68 in \cite{edwards1984galois}). In this article we want to understand how big or how small the subset of roots be so that it generates the splitting field. \smallskip

We are not interested in writing the splitting field $K_f$ as $K(\beta)$ for some primitive element $\beta\in K_f$. Our interest is writing this as $K_f=K(\beta_1, \beta_2, \dots, \beta_m)$ where $\{\beta_i\}_{i=1}^m$ is a subset of the roots of $f$ in $\bar{K}$. Such a subset is called a generating set of the splitting field of the polynomial $f$ over $K$. In this context it is therefore natural to study the subset of roots of $f$ belonging to the field generated by a single specific root of $f$. Inspired from the work of Perlis in \cite{perlis2004roots}, such subsets were defined as root clusters by the second author and Krithika in their recent work in \cite{krithika2023root}. So clusters form equivalence classes for the following equivalence relation among the roots of $f$ : $\alpha\sim \beta$ over $K$ if $K(\alpha)=K(\beta)$ (equality of fields and not mere isomorphism). For example all the primitive $n$-th roots of unity are equivalent to each other over $\mathbb{Q}$ but any two distinct cube roots of $2$ are not equivalent over $\mathbb{Q}$. The second author and Krithika also introduced the concept of cluster towers in \cite{krithika2023root}. In that work, they had posed a question about degree sequence of cluster tower being independent of the ordering of representatives of clusters of roots. This was answered by the first author and Bhagwat in their recent work on root clusters in \cite{Bhagwat_2025} by constructing a counterexample.


\smallskip
By a minimal generating set we mean a generating set for which no proper subset is a generating set. In this article, we prove some basic properties of minimal generating set of the splitting field of a polynomial over $K$. We then go on to establish Proposition \ref{equivalence} which connects the concepts of minimal generating sets and cluster towers and which aids in reformulating results on minimal generating sets as results on cluster towers. The notion of minimal cluster tower is also introduced. A word of caution: In the context of vector spaces we have a similar notion of minimal generating set. The Steinitz Exchange Lemma allows us to conclude that all such minimal generating sets have the same cardinality, consequently we have a well-defined notion of the dimension of a vector space. One of the aims of this article is to show that minimal generating sets defined for splitting fields of polynomials over $K$ in a natural way as given above does not possess the Steinitz Exchange property. In contrast we will be exhibiting two ``bases" with one of them of cardinality $2$ and other of given cardinality $k>2$. The precise result is established in Proposition \ref{two minimal sets different cardinality}.\smallskip

We establish the following result in this article. 

\begin{mainthm}
    
[Cardinalities of Minimal generating sets]\label{refer} \hfill

\begin{enumerate}

  \item  Let $n>2$ be an odd composite number. Fix $\zeta$ to be a primitive $n$-th root 
of unity in $\bar{\mathbb{Q}}$. Let $c$ be a positive rational number such that $f=x^n-c$ is
  an irreducible polynomial over $\mathbb{Q}$. Then the splitting field $\mathbb{Q}_f$ of $f$ over $\mathbb{Q}$ has minimal generating sets of cardinalities $2,3,\dots, \omega(n)$ (and minimal cluster towers of lengths $2,3,\dots, \omega(n)$) and these are the only possible cardinalities for minimal generating set for $f$ over $\mathbb{Q}$ (and these are the only possible lengths for minimal cluster towers for $f$). Here $\omega(n)$ denotes the number of distinct prime divisors of $n$.   \smallskip
  
\item  Let $K$ be a number field. For any $n>2$, there exist infinitely many degree $n$ irreducible polynomials over $K$ for which the splitting field has all its minimal generating sets of cardinality $k$ (and all cluster towers of length $k$) for the following values of $k$ : $(1)$ $k=n-1$ and $(2)$ $k=n-2$.\smallskip

\item There exists a degree $n$ irreducible polynomial over $\mathbb{Q}$ for which the splitting field has all its minimal generating sets of cardinality $k$ (and all cluster towers of length $k$) for the following values of $n$ and $k$:

\begin{enumerate}

\item (i) $n=12$ and $k=5$, (ii) $n=11$ and $k=4$.

\item $n=p+1$ where $p$ is an odd prime and $k=3$.

\item $n=q$ where $q>2$ is a prime power and $k=2$.
\end{enumerate}
\end{enumerate}
 
\end{mainthm}

Thus it is clear that our notion of a minimal generating set is different from the notion of a minimal strong base in \cite{mckay1997finding} (which is minimal system of roots in \cite{gomez2011sharply}) which we refer to as a shortest minimal generating set in this article i.e. a minimal generating set with least cardinality.\smallskip

Proposition \ref{deg seq thm} establishes interesting properties of length $k$ cluster tower associated with the cardinality $k$ minimal generating set that we constructed in proof of Proposition \ref{two minimal sets different cardinality}. As a consequence we get that even if we work with minimal generating set, the degree sequence still depends on the ordering of the elements in the set. Lemma \ref{iff for gen set} (1) establishes equivalent condition for a set to be a generating set of splitting field of $f$ over $\mathbb{Q}$ under a hypothesis similar to hypothesis in \MainTheorem{} (1). Under the same hypothesis, we have Lemma \ref{iff for gen set} (3) and (4) which facilitate the proof of \MainTheorem{} (1) and we also have Lemma \ref{iff for gen set} (2) where we establish an equivalent condition for a set to be shortest minimal generating set and count the total number of shortest minimal generating sets of splitting field of $f$ over $\mathbb{Q}$. Theorems \ref{mult transtvty} and \ref{sharply mult trans} establish how minimal generating sets and related concepts behave under multiple transitivity of the Galois group and as a consequence of these results we obtain \MainTheorem{} (2), (3).


  

  


    



\section{Minimal Generating Sets and Cluster Towers}\label{minimal section}

Let $K$ be a perfect field. We fix an algebraic closure $\bar{K}$ once and for all and work with irreducible polynomials over $K$ and finite extensions of $K$ contained in $\bar K$. Let $f \in K[t]$ be an irreducible polynomial of degree $n$. Since $K$ is perfect, it follows that $f$ has $n$ distinct roots in $\bar{K}$. We denote the set of roots of $f$ by $R=\{\alpha_1,\alpha_2,\dots,\alpha_n\}$. Let $K_f$ be the splitting field of $f$ over $K$ and let $G = {\rm {{\rm Gal}}}(K_f/K)$. 
We have the following notions related to root clusters of $f$ (For details see Sections 1 and 2 in \cite{krithika2023root} and Section 2.1 in \cite{Bhagwat_2025}): A cluster of $f$ is defined as the subset of roots of $f$ belonging to the field generated by a single root of $f$ over $K$. It turns out that all the clusters of $f$ have the same cardinality which is defined as the cluster size of $f$ over $K$ denoted as $r_K(f)$. The number of clusters of $f$ over $K$ is denoted as $s_K(f)$.

\subsection{Minimal generating sets of splitting field of a polynomial over a base field}  

\hfill\smallskip

For a set $B=\{\beta_1, \beta_2, \dots, \beta_m\}\subset \bar{K}$, let $K(B)$ denote the field $K(\beta_1, \beta_2, \dots, \beta_m)$.

  \begin{definition}
      Consider an irreducible polynomial $f$ over $K$. A set $B\subset R$ is called a minimal generating set of the splitting field of the polynomial $f$ over $K$ if the following hold.
      
      \begin{enumerate}
          \item $K(B)=K_f$  

          \item For any set $A\subsetneq B$, we have $K(A)\neq K_f$. 
      \end{enumerate}

We also refer to this set simply as a minimal generating set for $f$ over $K$.
      
      \end{definition}

\begin{proposition} \label{first prop} 
Let $B\subset R$ be a minimal generating set of the splitting field of the polynomial $f$ over $K$. Then we have the following
    \begin{enumerate}
        \item All the elements of $B$ lie in distinct root clusters of $f$. Thus $|B|\leq s_K(f)$. 

        \item $B\subset B'$ where $B'$ is a complete set of representatives of root clusters of $f$. 
    \end{enumerate}
\end{proposition}

\begin{proof}\hfill
    \begin{enumerate}
        \item If two elements of $B$, say $\beta_1$ and $\beta_2$ lie in the same root cluster, then $K(\beta_1)=K(\beta_2)$. Let $A=B\backslash\{\beta_1\}$. Then $A\subsetneq B$ with $K(A)=K_f$ which is a contradiction. Since $s_K(f)$ is the number of distinct clusters of $f$, we have $|B|\leq s_K(f)$.

        \item From part (1), we have that elements of $B$ are representatives of $|B|$ many clusters of $f$. We can choose representatives of remaining $s_K(f)-|B|$ many clusters and call that set $B_0$. Then $B'=B\cup B_0$ serves our purpose.
\end{enumerate}\end{proof}

The following lemma consists of some simple yet useful observations which are easy to verify.

\begin{lemma}
    \label{minimal equivalent def}
\hfill

\begin{enumerate}
\item 
Consider $B=\{\beta_1, \beta_2, \dots, \beta_m\}\subset R$. We have that $B$ is a minimal generating set of the splitting field of the polynomial $f$ over $K$, if and only if $K(B)=K_f$ and for every $1\leq i\leq m$, we have $\beta_i\not \in K(B\backslash \{\beta_i\})$. \smallskip

\item  Let $C\subset R$ be a generating set of the splitting field of polynomial $f$ over $K$ i.e. $K(C)=K_f$. Then there exists $B\subset C$ such that $B$ is a minimal generating set of the splitting field of $f$ over $K$. (In the context of vector spaces, we have that every spanning set contains a basis. We have a similar property holding true for minimal generating sets as well).
\smallskip

\item Let $f$ be an irreducible polynomial over $K$. Then there exists a minimal generating set $B$ of splitting field of $f$ over $K$.   
\end{enumerate}  
\end{lemma}

\begin{remark}
    For a minimal generating set $B$ of the splitting field of a polynomial $f$ over $K$, $|B|$ need not divide degree $n$ of $f$. Consider $f=x^3-2$ over $\mathbb{Q}$ and let $\omega$ be primitive $3$rd root of unity. Then $\{ \sqrt[3]{2}, \sqrt[3]{2}\omega\}$ is a minimal generating set for $f$ over $\mathbb{Q}$ but $2\nmid 3$.
\end{remark}


 \begin{proposition}
     
\label{two minimal sets different cardinality}
  Given an integer $k > 2$, there exist infinitely many irreducible polynomials $f$ over $\mathbb{Q}$ for which the splitting field $\mathbb{Q}_f$ for each of the polynomial $f$ over $\mathbb{Q}$ has two minimal generating sets of different cardinalities: one of cardinality 2 and another of cardinality k. In fact we can obtain infinitely many such polynomials of the same degree as well.
 \end{proposition}

\begin{proof}

Let $n=p_1 p_2 \cdots p_k$ where $p_i$'s are distinct odd primes. Fix $\zeta$ to be a primitive $n$-th root 
of unity in $\bar{\mathbb{Q}}$. Let $c$ be a positive rational number such that $f=x^n-c$ is
  an irreducible polynomial over $\mathbb{Q}$. Let $a = c^{1/n}$ be the positive real root. Thus $\mathbb{Q}_f=\mathbb{Q}(a,\zeta)$. Clearly $B_1=\{a,a\zeta\}$ is a minimal generating set for $f$ over $\mathbb{Q}$ with cardinality 2. We claim that $B_2=\{a\zeta^{n/p_i}\}_{i=1}^{k}$ is a minimal generating set for $f$ over $\mathbb{Q}$ with cardinality $k$.\smallskip

  Since $n$ is odd, by Proposition 1 in \cite{jacobson1990galois} and Theorem A in \cite{jacobson1990galois}, $\mathbb{Q}(a)\cap \mathbb{Q}(\zeta)=\mathbb{Q}$ and $Gal(\mathbb{Q}_f/\mathbb{Q})\cong\mathbb{Z}/n \mathbb{Z} \rtimes (\mathbb{Z}/n \mathbb{Z})^{\times}$. Identifying $G={\rm Gal}(\mathbb{Q}_f/\mathbb{Q})$ with $\mathbb{Z}/n \mathbb{Z} \rtimes (\mathbb{Z}/n \mathbb{Z})^{\times}$ we have $H=Gal(\mathbb{Q}_f/\mathbb{Q}(a))= \{0\}\times (\mathbb{Z}/n \mathbb{Z})^{\times}\subseteq G$. Now $G$ has the semidirect product group law \[(\alpha,u)\cdot (\beta,v)=(\alpha + u\cdot \beta, uv)\] 
as in \cite{jacobson1990galois} where $u\cdot\beta$ is usual multiplication $u\beta$ in the ring $\mathbb{Z}/n \mathbb{Z}$. The action of $G$ on the roots of $f$ is given by $(\alpha, u).(a\zeta^j)=a \zeta^{(\alpha+ uj)}$ for any $1\leq j\leq n$. 
Thus we observe that for each $j$, $H_{j}=\{(j-vj,v)\ |\ v\in (\mathbb{Z}/n \mathbb{Z})^{\times}\}$ is the subgroup of $G$ fixing the field $\mathbb{Q}(a\zeta^{j})$. \smallskip

Consider $(\alpha,u)\in \cap_{i=1}^k\ H_{n/p_i}$. Thus $\alpha=(1-u)n/p_i$ for all $1\leq i\leq k$. Hence for each $i$, $p_i| \alpha$. Therefore $n|\alpha$ which means $\alpha=0$. Thus for each $i$, $p_i|(1-u)$. Therefore $n|(1-u)$ which means $u=1$. Thus we have shown $\cap_{i=1}^k\ H_{n/p_i}=1$ which means $K(B_2)=\mathbb{Q}_f$.\smallskip

Now for any $1\leq l\leq k$, consider $(\alpha,u)\in \cap_{i\neq l}\ H_{n/p_i}$. Thus $\alpha=(1-u)n/p_i$ for all $i\neq l$. Hence for each $1\leq i\leq k$, $p_i| \alpha$. Therefore $\alpha=0$. Thus for each $i\neq l$, $p_i|(1-u)$. Therefore $(n/p_l) |(1-u)$. We claim that there exists an integer $1\leq z\leq (p_l -1)$ such that $(n/p_l)z+1\in (\mathbb{Z}/n \mathbb{Z})^{\times}$. \smallskip

It is enough to show that there is a $1\leq z\leq (p_l-1)$ such that $p_l\nmid ((n/p_l)z+1)$. Assume on the contrary that $p_l | ((n/p_l)z+1)$ for all $1\leq z\leq (p_l-1)$. Thus $p_l$ divides the sum $\Sigma_{z=1}^{(p_l-1)}\ ((n/p_l)z+1)$ $=(n/p_l) (\Sigma_{z=1}^{(p_l-1)}\ z) + (p_l-1)=(n/p_l)(p_l)(p_l-1)/2 +(p_l-1)$ which is a contradiction. \smallskip

Hence we have a nontrivial element $(0, (n/p_l)z+1)\in \cap_{i\neq l}\ H_{n/p_i}$ for that choice of $z$. Thus for any $1\leq l\leq k$, we have $K(B_2\backslash \{a\zeta^{n/p_l}\})\neq \mathbb{Q}_f$. Hence by Lemma \ref{minimal equivalent def} (1), we are done. We get infinitely many polynomials with the property since we have infinite choices for $n=p_1 p_2 \cdots p_k$ and $c$. By fixing $n$ and varying $c$ over primes, we get by Eisenstein criterion, infinitely many such polynomials of same degree. \end{proof}

\begin{definition}
    $B\subset R$ is said to be shortest minimal generating set of the splitting field of the polynomial $f$ over $K$ if $B$ is minimal generating set with least possible cardinality. We similarly define longest minimal generating set as one having greatest possible cardinality.
\end{definition}

\begin{remark}
   Note that our notion of a shortest minimal generating set is same as the notion of a minimal strong base in \cite{mckay1997finding} (which is same as minimal system of roots in \cite{gomez2011sharply}).
\end{remark}

\subsection{Cluster tower of a polynomial over a base field and minimal generating sets}

In \cite{krithika2023root}, the notion of cluster towers is introduced. See Section 5.2 in \cite{Bhagwat_2025} for a group-theoretic formulation.\smallskip

 Let $f$ be an irreducible polynomial over $K$. Consider a complete set of representatives of clusters of roots of $f$ in $\bar{K}$. Let $(\beta_1, \beta_2,\ldots, \beta_s)$ be an ordering of this set where $s=s_K(f)$. Now consider the following cluster tower of fields
terminating at the splitting field $K_f$. 
\[ K \subseteq K(\beta_1) \subseteq K(\beta_1, \beta_2) \subseteq \dots \subseteq K(\beta_1, \beta_2, \ldots , \beta_s) = K_f.\]

The length of tower was defined in \cite{krithika2023root} as the number of distinct fields in the tower so that in the case when the Galois group of $f$ over $K$ is $\mathfrak S_n$ for $n > 2$, we have $n=$ number of clusters$=$ length of tower. However in this article we redefine the length of tower as the number of strict inclusions in the tower (in order to align this notion with standard conventions of lengths regarding chains of vector subspaces or prime ideals in the context of commutative rings and in order to ensure consistency in the results). Clearly length of tower $\leq s$. \smallskip

We also have the notion of degree sequence of the tower which is the ordered sequence of degrees of the distinct extensions over $K$ in the tower i.e. $([L_1:K], [L_2 :K],\dots,[K_f:K] )$ if the distinct fields in the tower are $K \subsetneq L_1 \subsetneq L_2 \subsetneq \dots \subsetneq  K_f$. Example 5.1.3 in \cite{Bhagwat_2025} demonstrates that both the degree sequence and length of tower are dependent on the ordering of the $\beta_i$’s.\smallskip

The following result connects the concepts of minimal generating sets and cluster towers.

\begin{proposition}

\label{equivalence}

Consider $B=\{\beta_1, \beta_2, \dots, \beta_m\}\subset R$.
    \begin{enumerate}
        \item $B$ is a minimal generating set of the splitting field of the polynomial $f$ over $K$, if and only if for every permutation $(i_1,i_2,\dots, i_m)$ of $(1,2,\dots, m)$, 
 \[  K \subseteq K(\beta_{i_1}) \subseteq K(\beta_{i_1}, \beta_{i_2}) \subseteq \dots \subseteq K(\beta_{i_1}, \beta_{i_2}, \ldots , \beta_{i_m})\] is a cluster tower for $f$ of length $m$.\smallskip

 In this case we refer any such cluster tower as a minimal cluster tower.\smallskip

 \item If the fields in two minimal cluster towers are equal at each step then we say that the two towers are equal, otherwise they are distinct. In the above case, all the $m!$ many minimal cluster towers are distinct. In other words, $B$ has $m!$ many minimal cluster towers associated with it.\smallskip

 \item   $B$ is a shortest minimal generating set of the splitting field of the polynomial $f$ over $K$, if and only if 
 \[  K \subseteq K(\beta_{1}) \subseteq K(\beta_{1}, \beta_{2}) \subseteq \dots \subseteq K(\beta_{1}, \beta_{2}, \ldots , \beta_{m})\] is a cluster tower for $f$ of length $m$ which is the least possible length for any tower for $f$.\smallskip

 In this case we refer such a cluster tower as a shortest minimal cluster tower.\end{enumerate}
\end{proposition}

\begin{proof} \hfill
\begin{enumerate}
    \item Suppose $B$ is a minimal generating set. Consider any permutation $(i_1,i_2,\dots, i_m)$ of $(1,2,\dots, m)$. Clearly the following is a cluster tower terminating at the splitting field  \[ K \subseteq K(\beta_{i_1}) \subseteq K(\beta_{i_1}, \beta_{i_2}) \subseteq \dots \subseteq K(\beta_{i_1}, \beta_{i_2}, \ldots , \beta_{i_m})=K_f\] since $K(B)=K_f$. If the length of this cluster tower is $<m$, then there exists $2\leq j\leq m$ such that $\beta_{i_j}\in K(\beta_{i_1},\beta_{i_2},\dots,\beta_{i_{j-1}})$. Let $A=B\backslash \{\beta_{i_j}\}$. We have $A\subsetneq B$ with $K(A)=K_f$ which contradicts $B$ being a minimal generating set. Hence the length of this cluster tower is $m$. Conversely, suppose $B$ is not a minimal generating set.\smallskip

    Case 1: $K(B)\neq K_f$. Then for any permutation, this tower is not a cluster tower since it doesn't terminate at the splitting field.

    Case 2: $K(B)= K_f$ and there exists $A\subsetneq B$ such that $K(A)=K_f$. Consider a permutation $(i_1,i_2,\dots, i_m)$ of $(1,2,\dots, m)$ such that $A=\{\beta_{i_1}, \beta_{i_2}, \dots, \beta_{i_l}\}$ where $l<m$. Then length of this cluster tower \[ K \subseteq K(\beta_{i_1}) \subseteq K(\beta_{i_1}, \beta_{i_2}) \subseteq \dots \subseteq K(\beta_{i_1}, \beta_{i_2}, \ldots , \beta_{i_m})=K_f\] is $\leq l< m$.
\smallskip

\item Suppose $(i_1,i_2,\dots,i_m)$ and $(j_1,j_2,\dots, j_m)$ are two distinct permutations of $(1,2,\dots,m)$. Consider the smallest $1\leq l \leq m-1$ such that $i_l\neq j_l$. We claim that $K(\beta_{i_1}, \beta_{i_2}, \ldots , \beta_{i_l})\neq K(\beta_{j_1}, \beta_{j_2}, \ldots , \beta_{j_l})$. Assume on the contrary, $K(\beta_{i_1}, \beta_{i_2}, \ldots , \beta_{i_l})=K(\beta_{j_1}, \beta_{j_2}, \ldots , \beta_{j_l})$. Since $i_k=j_k$ for all $1\leq k\leq l-1$. Thus $j_l\neq i_k$ for any $1\leq k\leq l$. This implies $\beta_{j_l}\in K(B\backslash \{\beta_{j_l}\})$ which by Lemma \ref{minimal equivalent def} (1), is a contradiction to $B$ being a minimal generating set. Hence the cluster towers corresponding to the two permutations are distinct.

\smallskip

  \item  Suppose \[  K \subseteq K(\beta_{1}) \subseteq K(\beta_{1}, \beta_{2}) \subseteq \dots \subseteq K(\beta_{1}, \beta_{2}, \ldots , \beta_{m})\] is a cluster tower for $f$ of length $m$ which is the least possible length for any tower for $f$. Thus for every permutation $(i_1,i_2,\dots, i_m)$ of $(1,2,\dots, m)$, 
 \[  K \subseteq K(\beta_{i_1}) \subseteq K(\beta_{i_1}, \beta_{i_2}) \subseteq \dots \subseteq K(\beta_{i_1}, \beta_{i_2}, \ldots , \beta_{i_m})\] is a cluster tower for $f$ of length $m$. By part (1), $B$ is a minimal generating set. If $B'$ is any other minimal generating set, then by part (1), the length of the corresponding cluster towers will be $|B'|$, since $|B|$ is the least possible length for any cluster tower for $f$. Hence $|B|\leq |B'|$. Hence $B$ is a shortest minimal generating set. Conversely, suppose $B$ is a shortest minimal generating set. Then by part (1), since $B$ is a minimal generating set, \[  K \subseteq K(\beta_{1}) \subseteq K(\beta_{1}, \beta_{2}) \subseteq \dots \subseteq K(\beta_{1}, \beta_{2}, \ldots , \beta_{m})\] is a cluster tower for $f$ of length $m$. Suppose $m$ is not the least possible length for cluster towers for $f$. Thus there is a set $C=\{\gamma_1,\gamma_2,\dots, \gamma_{m'}\}\subset R$ with $m'< m$ such that \[  K \subseteq K(\gamma_{1}) \subseteq K(\gamma_{1}, \gamma_{2}) \subseteq \dots \subseteq K(\gamma_{1}, \gamma_{2}, \ldots , \gamma_{m'})\] is a cluster tower for $f$ of length $m'<m$. Thus $K(C)=K_f$. Hence by Lemma \ref{minimal equivalent def} (2), there exists $B'\subset C$ such that $B'$ is a minimal generating set of the splitting field of $f$ over $K$. Thus $|B'|\leq |C| <m$. This gives a contradiction to $B$ being a shortest minimal generating set.\end{enumerate}\end{proof}


\begin{corollary}
    Let $f$ be an irreducible polynomial over $K$ and $\alpha$ be a root of $f$ in $\bar{K}$. Then $f$ has a unique cluster tower if and only if $K(\alpha)/K$ is Galois.\smallskip

     In this case, the length of the tower is $1$ and degree sequence is $(deg(f))$ and there are $deg(f)$ many minimal generating sets which are singletons and hence also shortest minimal generating sets.
\end{corollary}

\begin{proposition}
    
 \label{deg seq thm}
Consider the minimal generating set $B_2=\{a\zeta^{n/p_i}\}_{i=1}^{k}$ with cardinality $k$ as in proof of Proposition \ref{two minimal sets different cardinality}. Consider the following length $k$ cluster tower associated with it. \[  \mathbb{Q} \subseteq L_1 \subseteq L_2 \subseteq \dots \subseteq L_k =\mathbb{Q}_f\] where $L_i=\mathbb{Q}(a\zeta^{n/p_1},a\zeta^{n/p_2},\dots, a\zeta^{n/p_i})$ for $1\leq i \leq k$. Then we have the following.

\smallskip

\begin{enumerate}
    \item $L_i=\mathbb{Q}(a, \zeta^{n/(p_1 p_2\cdots p_i)})$ for all $i\geq 2$.\smallskip
    
    \item The degree sequence of the cluster tower is\\
    $(n, n\phi(p_1p_2), n\phi(p_1p_2p_3),\dots, n \phi(p_1p_2\cdots p_{k-1}),n \phi(n))$.\smallskip


    \item The root capacity of $L_i$ with respect to $a$ with base field $\mathbb{Q}$ i.e. the number of roots of minimal polynomial of $a$ over $\mathbb{Q}$ that are contained in $L_i$ is $\rho_{\mathbb{Q}} (L_i,a)=p_1p_2\cdots p_i$ for $i\geq 2$ and $\rho_{\mathbb{Q}} (L_1,a)=1$.
    (Refer Section 6.2 in \cite{Bhagwat_2025} for more on the concept of root capacity).

\end{enumerate}

\end{proposition}

\begin{proof}\hfill
    \begin{enumerate}
    
 \item We have $i\geq 2$. 
 By argument in proof of Lemma \ref{iff for gen set} (1) below, $L_i=\mathbb{Q}(a\zeta^{n/p_1},\zeta^l)$ where $$l=gcd(n, (n/p_2 - n/p_1), (n/p_3 - n/p_1),\dots, (n/p_i - n/p_1) )=n/(p_1p_2\cdots p_i).$$ Since $l|(n/p_1)$, we have $L_i=\mathbb{Q}(a,\zeta^l)$.
 

 \smallskip
 
\item Now $\mathbb{Q}(a)\cap \mathbb{Q}(\zeta)=\mathbb{Q}$. Thus $\mathbb{Q}(a)\cap \mathbb{Q}(\zeta^{n/(p_1 p_2\cdots p_i)})=\mathbb{Q}$. Hence $[\mathbb{Q}(a, \zeta^{n/(p_1 p_2\cdots p_i)}): \mathbb{Q}]=[\mathbb{Q}(a):\mathbb{Q}][\mathbb{Q}(\zeta^{n/(p_1 p_2\cdots p_i)}):\mathbb{Q}]$. Thus for $i\geq 2$, $[L_i : \mathbb{Q}]=n\phi(p_1p_2\cdots p_i)$.

\smallskip


\item Since the field $\mathbb{Q}(a)$ contains a single root $a$ of $f=x^n-c$ as $n$ is odd. Thus the field $\mathbb{Q}(a\zeta^{n/p_1})$ which is isomorphic to $\mathbb{Q}(a)$ also contains a single root $a\zeta^{n/p_1}$ of $f$. Hence $\rho_{\mathbb{Q}}(L_1, a)=1$. Also since $L_k=L_f$. Hence $\rho_{\mathbb{Q}}(L_k,a)=n=p_1p_2\cdots p_k$. 

\smallskip

Now let $2\leq i\leq k-1$. Clearly $p_1p_2\cdots p_i$ many roots $\{a\zeta^{yn/(p_1p_2\cdots p_i)}\}_{y=0}^{(p_1p_2 \cdots p_i -1)}$ of $f$ lie in $L_i$. We will show that these are the only roots that lie in $L_i$. Suppose some other root $a\zeta^j$ of $f$ also lies in $L_i$. Hence there exists an $i+1\leq l\leq k$ such that $p_l\nmid j$. Since $L_i=\mathbb{Q}(a\zeta^{n/p_1},a\zeta^{n/p_2},\dots, a\zeta^{n/p_i})$. Thus subgroup of $G$ fixing $L_i$ is $\cap_{w=1}^i H_{n/p_w}$. Since $\mathbb{Q}(a\zeta^j)\subset L_i$. Thus $\cap_{w=1}^i H_{n/p_w}\subset H_j$. Hence the intersection of $k-1$ subgroups $\cap_{w\neq l}\ H_{n/p_w}\subset H_j$. Now in last paragraph of proof of Proposition \ref{two minimal sets different cardinality} we got a nontrivial element $(0, (n/p_l)z+1)\in \cap_{w\neq l}\ H_{n/p_w}$ for some $1\leq z\leq p_l-1$. Thus $(0, (n/p_l)z+1)\in H_j$. Hence $n| ((n/p_l)zj)$. But since $p_l\nmid ((n/p_l)zj)$, we arrive at a contradiction.\end{enumerate}\end{proof}


\begin{remark}
      In Proposition \ref{deg seq thm}, if instead we consider the cluster tower corresponding to the permutation $(1,3,2,4,\dots,k)$ of $(1,2,3,4,\dots,k)$, we get the degree sequence as\\
      $(n, n\phi(p_1p_3), n\phi(p_1p_2p_3),\dots, n \phi(p_1p_2\cdots p_{k-1}),n \phi(n))$. Now $\phi(p_1p_3)\neq \phi(p_1p_2)$. Hence the degree sequences for both the cluster towers are distinct. This demonstrates that even if we work with minimal generating set, the degree sequence still depends on the ordering of the elements in the set.
\end{remark}

\subsection{A special family of polynomials}

Inspired by the proof of Proposition \ref{two minimal sets different cardinality}, in the subsequent discussion we will study a special family of polynomials more closely.

\begin{lemma}\label{iff for gen set}
     Let $n>2$ be odd. Fix $\zeta$ to be a primitive $n$-th root 
of unity in $\bar{\mathbb{Q}}$. Let $c$ be a positive rational number such that $f=x^n-c$ is
  an irreducible polynomial over $\mathbb{Q}$. Let $a = c^{1/n}$ be the positive real root. 
  
\begin{enumerate}
\item  $C=\{a\zeta^{b_i}\}_{i=1}^m$ is a generating set for $f$ over $\mathbb{Q}$ if and only if $$gcd(n, b_2-b_1, b_3-b_1,\dots,b_m-b_1)=1.$$

\item  
 $B=\{a\zeta^k,a\zeta^l \}$ for $0\leq k < l \leq n-1$ is a shortest minimal generating set for $f$ over $\mathbb{Q}$ if and only if $gcd(l-k, n)=1$. Thus there are $n\phi(n)/2$ many shortest minimal generating sets for $f$ over $\mathbb{Q}$. 

\smallskip

\item  Suppose for $m>2$, $C=\{a\zeta^{b_i}\}_{i=1}^m$ is a minimal generating set for $f$ over $\mathbb{Q}$. Then $m\leq \omega(n)$ where $\omega(n)$ is the number of distinct prime divisors of $n$. \smallskip

\item If for some $m\geq 2$ we have $n=a_1a_2\cdots a_m$ where $a_i$'s are pairwise coprime numbers greater than 1, then $B=\{a\zeta^{n/a_i}\}_{i=1}^m$ is a minimal generating set for $f$ over $\mathbb{Q}$. 

\end{enumerate} \end{lemma}

\begin{proof}
\hfill

\begin{enumerate}
\item Let $K=\mathbb{Q}(a\zeta^{b_1},a\zeta^{b_2},\dots, a\zeta^{b_m})$. Clearly $K=\mathbb{Q}(a\zeta^{b_1}, \zeta^{b_2-b_1},\zeta^{b_3-b_1},\dots,\zeta^{b_m-b_1})$. Also we have $K\subset \mathbb{Q}(a\zeta^{b_1}, \zeta^l)$ where $l=gcd(n, b_2-b_1, b_3-b_1,\dots,b_m-b_1)$. We have integers $y,x_2,x_3,\dots,x_m$ such that $y n + \Sigma_{j=2}^m x_j (b_2 - b_1)= l$. Hence $K=\mathbb{Q}(a\zeta^{b_1},\zeta^l)$. We know $\mathbb{Q}_f=\mathbb{Q}(a,\zeta)$. Suppose $l=1$. Then clearly $K=\mathbb{Q}_f$.\smallskip
   
   Conversely suppose $K=\mathbb{Q}_f$. Since $n$ is odd, by results in \cite{jacobson1990galois}, we have $\mathbb{Q}(a)\cap \mathbb{Q}(\zeta)=\mathbb{Q}$ (which is equivalent to $\mathbb{Q}(a)$ and $\mathbb{Q}(\zeta)$ being linearly disjoint over $\mathbb{Q}$ by Example 20.6 in \cite{morandi2012field}). Therefore $[\mathbb{Q}_f:\mathbb{Q}]=n\phi(n)$ where $\phi$ is the Euler totient function. Now since $\mathbb{Q}(a\zeta^{b_1}, \zeta)=\mathbb{Q}_f$. Thus $\mathbb{Q}(a\zeta^{b_1}) \cap \mathbb{Q}(\zeta)=\mathbb{Q}$. Hence $\mathbb{Q}(a\zeta^{b_1})\cap \mathbb{Q}(\zeta^l)=\mathbb{Q}$. Also since $l|n$, we have $[\mathbb{Q}(\zeta^l):\mathbb{Q}]=\phi(n/l)$. Thus $[K:\mathbb{Q}]=n\phi(n/l)$. Since $K=\mathbb{Q}_f$. Thus $\phi(n)=\phi(n/l)$. Since $n$ is odd, $n=n/l$ i.e. $l=1$.\smallskip

\item Clearly $\{a\zeta^k\}$ for any $0\leq k\leq n-1$ is not a generating set for $f$ over $\mathbb{Q}$. Thus $B=\{a\zeta^k,a\zeta^l \}$ for $0\leq k < l \leq n-1$ is a shortest minimal generating set for $f$ over $\mathbb{Q}$ if and only if $B=\{a\zeta^k,a\zeta^l \}$ is a generating set for $f$ over $\mathbb{Q}$ which is equivalent to $gcd(l-k, n)=1$ by part (1).\smallskip

We will count the number of pairs $(k,l)$ such that $0\leq k < l \leq n-1$ and $gcd(l-k, n)=1$. Since $n>2$, $\phi(n)$ is even. Let $e=\phi(n)$ and let the $e$ many elements in $(\mathbb{Z}/n\mathbb{Z})^{\times}$ be \\
$x_1=1<x_2<\dots<x_e=(n-1)$. Now for a given $k$, we must have $l=k+x_i$ for some $1\leq i \leq e$ so that the gcd condition is satisfied. When $k=0$, then $l$ has $e$ choices. Consider $0\leq y\leq (e-2)$. 
For $(n-x_{e-y}) \leq k \leq (n-1-x_{e-y-1})$ we have $e-y-1$ choices for $l$. Hence the total number of pairs is $e + \Sigma_{y=0}^{e-2}\ (x_{e-y}-x_{e-y-1})(e-y-1)=ne-\Sigma_{i=1}^e\ x_i$. Observing that $x_{e-i+1}=n-x_{i}$ for all $1\leq i\leq e/2$, we have $\Sigma_{i=1}^e\ x_i=ne/2$ and thus we are done.
  \smallskip

\item Since $C$ is a generating set. Thus from part (1) we have $$gcd(n, b_2-b_1, b_3-b_1,\dots,b_m-b_1)=1.$$ Since $C$ is minimal generating set, for any $1\leq j\leq m$ we have that $C\backslash\{a\zeta^{b_j}\}$ is not a generating set. Thus for $2\leq j\leq m$ we get by part (1) that $$gcd(n, b_2-b_1,\dots, b_{j-1}-b_1,b_{j+1}-b_1, \dots,b_m-b_1)\neq 1.$$ Thus for any $2\leq j\leq m$, we have that there exists a prime $p_j|n$ such that $p_j|(b_l-b_1)$ for all $l\neq j$ and $p_j\nmid (b_j-b_1)$. Clearly these $m-1$ many primes are distinct. Observe that since $C$ is a generating set, we also have $$gcd(n, b_1-b_2, b_3-b_2,\dots,b_m-b_2)=1.$$ Since $C\backslash\{a\zeta^{b_1}\}$ is not a generating set, we have $$gcd(n, b_3-b_2, b_4-b_2,\dots,b_m-b_2)\neq 1.$$ Hence we have a prime $p_1|n$ such that $p_1|(b_l-b_2)$ for all $l\neq 1$ and $p_1\nmid (b_1-b_2)$. Now for $3\leq j\leq m$ we have $p_j|(b_2-b_1)$ and thus $p_1\neq p_j$. Now suppose $p_1=p_2$. We have $p_2|(b_3-b_1)$ and $p_1|(b_3-b_2)$. Thus $p_1|((b_3-b_2)-(b_3-b_1))$ that is $p_1 |(b_1-b_2)$ which is a contradiction. Thus $p_1\neq p_2$. Thus $n$ has atleast $m$ many distinct prime divisors.

    \smallskip
    
\item  Let $l=gcd(n, n/a_2-n/a_1, n/a_3-n/a_1,\dots,n/a_m-n/a_1)$. Now $l|a_1a_2\cdots a_m$. Now for any $2\leq j\leq m$, we have $p_j\nmid (n/a_j-n/a_1)$ where $p_j$ is any prime dividing $a_j$. Thus $p_j\nmid l$. Similarly we have $p_1\nmid (n/a_2-n/a_1)$ for any prime $p_1|a_1$. Thus $p_1\nmid l$. Hence $l=1$. Thus $B$ is a generating set for $f$ over $\mathbb{Q}$ by part (1).\smallskip

 We will show that $B\backslash \{a\zeta^{n/a_m}\}$ is not a generating set. Similarly one can show that $B\backslash \{a\zeta^{n/a_j}\}$ is not a generating set for $1\leq j\leq m-1$ and thus $B$ is a minimal generating set. Let $d=gcd(n, n/a_2-n/a_1, n/a_3-n/a_1,\dots,n/a_{m-1}-n/a_1)$. Now $d|n$. Clearly $a_m|(n/a_j-n/a_1)$ for every $2\leq j\leq m-1$. Thus $a_m|d$. By similar arguments as in first paragraph we have for any $1\leq j\leq m-1$ that $p_j\nmid d$ for any prime $p_j|a_j$. Thus $d=a_m\neq 1$. Thus by part (1) we are done.\end{enumerate}\end{proof}



 \begin{remark} \label{prime}
\hfill

 \begin{enumerate}
     \item 
 
If $n$ is an odd prime $p$ in Lemma \ref{iff for gen set} (2), then the $p\phi(p)/2 = {p\choose 2}$ many shortest minimal generating sets for $f$ over $\mathbb{Q}$, $\{a\zeta^k,a\zeta^l \}$ for $0\leq k < l \leq p-1$ are all the possible minimal generating sets for $f$ over $\mathbb{Q}$. \smallskip

\item Consider the case in Lemma \ref{iff for gen set} (2). Then by Proposition \ref{equivalence} (3), there are totally $n$ many distinct shortest minimal cluster towers of $f$. For each tower, the length is $2$ and its degree sequence is $(n,n\phi(n))$. If $n$ is an odd prime $p$, then there are totally $p$ many distinct cluster towers of $f$. 

\smallskip

\item   Lemma \ref{iff for gen set} (3) is not true for $m=2$ as for $n=p$ prime we have $\omega(n)=1<2$ and any minimal generating set has cardinality $2$.\smallskip

\item  Consider the case in proof of Proposition \ref{two minimal sets different cardinality}. The minimal generating set of cardinality $k=\omega(n)$ is a longest minimal generating set of the splitting field of $f$ over $\mathbb{Q}$. \smallskip

\item  One can observe that Lemma \ref{iff for gen set} (4) gives an alternate proof for $B_2$ being a minimal generating set in Proposition \ref{two minimal sets different cardinality}.

\end{enumerate}
 \end{remark}




 



\begin{proof}[Proof of \MainTheorem{} (1)]
    By Lemma \ref{iff for gen set} (3), the cardinality of any minimal generating set is bounded by $\omega(n)$. Let $k=\omega(n)$. So we have $n=p_1^{c_1}p_2^{c_2}\cdots p_k^{c_k}$. For any $2\leq l\leq k$, let $a_1=p_1^{c_1}, a_2=p_2^{c_2}, \dots , a_{l-1}=p_{l-1}^{c_{l-1}}$ and $a_l=p_l^{c_l}p_{l+1}^{c_{l+1}}\dots p_k^{c_k}$. Hence by Lemma \ref{iff for gen set} (4), we get a minimal generating set for $f$ over $\mathbb{Q}$ of cardinality $l$. The assertions about cluster towers follow from Proposition \ref{equivalence}.\end{proof}

\begin{remark}
In the context of vector spaces, we have that every linearly independent set can be extended to a basis. We don't have a similar property holding true for minimal generating sets. Let $f$ be an $n$ degree irreducible polynomial over a perfect field $K$ and let $R$ be the set of its roots. Let $D=\{\delta_1,\delta_2,\dots, \delta_{m'}\} \subset R$ such that $\delta_i\not \in K(D\backslash \{\delta_i\})$ for all $1\leq i \leq m'$. Then it is not necessary that there exists a $B\supset D$ such that $B$ is a minimal generating set of the splitting field of $f$ over $K$.\smallskip

Consider the case in \MainTheorem{} (1) for $n=105=3\cdot 5\cdot7$ with $D=\{a,a\zeta^{15}\}$. Clearly $a\zeta^{15}\not \in K(a)$ and $a\not \in K(a\zeta^{15})$. Since $gcd(15-0,105)=15\neq 1$, $D$ is itself not a minimal generating set. Suppose $B\supset D$ such that $B$ is a minimal generating set. Then $|B|=3$. Let $B=D\cup \{a\zeta^j \}$ where $j\neq 0, 15$. We should necessarily have the following: $gcd(j-0,105)\neq 1$, $gcd(j-15,105)\neq 1$ and $gcd(105, 15-0, j-0)=1$.\smallskip

Now $1=gcd(105,15,j)=gcd(15,j)$. Thus $3,5\nmid j$. Also since $gcd(j,105)\neq 1$, we have $j=7l$ for some $1\leq l<15$ and $3,5\nmid l$. Now we have $gcd(j-15, 105)\neq 1$. But $gcd(7l-15,105)=1$ since $3,5\nmid l$. This gives a contradiction.
\end{remark}

\section{Multiple transitivity and Minimal generating sets}

We begin by recalling the following important and useful result which is the final proposition in Perlis \cite{perlisroots} which establishes that every faithful transitive permutation representation
of an abstract Galois group occurs concretely as the root permutation group
of a polynomial.

\begin{proposition}\label{perm}
 Let $G$ be a transitive subgroup of ${\mathfrak S}_n$ for some $n$. If there exists a finite Galois extension of a field $K$ with Galois group isomorphic to $G$, then there exists an irreducible polynomial $f$ over $K$ of degree $n$
and a labelling of the roots of $f$ so that the Galois group of $f$, viewed as a group permuting the roots of $f$, is precisely $G$.

\end{proposition}


\begin{proof}[Proof of \MainTheorem{} (2) for $k=n-1$]

    By results in \cite{volklein1996groups} on hilbertian fields, we have ${\mathfrak S}_{n}$ to be realizable as a Galois group for infinitely many pairwise linearly disjoint Galois extensions over $K$. Thus by Proposition \ref{perm}, there exist infinitely many irreducible polynomials $f$ over $K$ of degree $n$ with Galois group ${\mathfrak S}_{n}$. Let the set of roots of $f$ be $R=\{\alpha_i\}_{i=1}^{n}\subset \bar{K}$. Let $B_i=R\backslash \{\alpha_i\}$ for $1\leq i\leq n$. Then each $B_i$ is a minimal generating set with cardinality $n-1$. \end{proof}
    
\begin{remark}
    In the above proof, the $n$ many $B_i$'s are all the possible minimal generating sets of the splitting field of $f$ over $K$. By a similar argument as in proof of Proposition \ref{equivalence} (2), we can show that, there are totally $n!/2$ many distinct cluster towers of $f$.
\end{remark}

We now recall an important notion.

\begin{definition}
    Consider an irreducible polynomial $f$ of degree $n>2$ over a perfect field $K$ with the Galois group $G=Gal(K_f/K)$. We say that $G$ is $k$-fold transitive (i.e. the action of $G$ on the set of all roots of $f$ in $K_f$ is $k$-fold transitive) if for any pair of cardinality-$k$ ordered subsets of roots of $f$, say $\{\alpha_1,\alpha_2,\dots, \alpha_k\}$ and $\{\beta_1,\beta_2,\dots, \beta_k\}$, we have a $\sigma\in G$ such that $\sigma(\alpha_i)=\beta_i$ for all $1\leq i\leq k$. We often say $k$-transitive to mean $k$-fold transitive. \smallskip

    Further if such a $\sigma$ is unique for each pair of cardinality-$k$ ordered subsets, then we say that $G$ is sharply $k$-transitive (i.e. the action of $G$ on the set of all roots of $f$ in $K_f$ is sharply $k$-transitive).
\end{definition}

The following result demonstrates how minimal generating sets and the related concepts behave under the assumption of multiple transitivity of the Galois group.

\begin{theorem}
    
\label{mult transtvty}

Consider an irreducible polynomial $f$ of degree $n>2$ over a perfect field $K$ with the Galois group $G=Gal(K_f/K)$. Suppose $k<n$. Then $G$ is $k$-transitive if and only if we have all of the following.

\smallskip

\begin{enumerate}
    \item Every minimal generating set for $f$ over $K$ has cardinality $\geq k$. \smallskip

 \item Every subset $D$ of roots of $f$ of cardinality $\leq k$ satisfies $\delta\not \in K(D\backslash\{\delta\})$ for all $\delta\in D$. Also every such $D$ is contained in some minimal generating set $B$ for $f$ over $K$.

 \smallskip

 \item If $k\geq 2$, then $r_K(f)=1$ and $s_K(f)=n$.\smallskip

 \item The length of any cluster tower of $f$ is $\geq k$ and the first $k$ terms in the degree sequence are $\ ^nP_1,\ ^nP_2 , \ldots,\ ^nP_{k}$\ . Thus there are $\geq\ ^nP_{k-1}$ many distinct cluster towers.

\end{enumerate}
    
\end{theorem}

\begin{proof}We first assume that $G$ is $k$-transitive and prove parts (1)-(4).

\begin{enumerate}
 \item The assertion is trivial for $k=1$. So assume $k>1$. Suppose on the contrary we have a minimal generating set for $f$ over $K$ say $B=\{\beta_1,\beta_2,\dots, \beta_m\}$ such that $m<k$. Since $k<n$, we have atleast two distinct roots of $f$ say $\alpha$ and $\alpha'$ which do not lie in $B$. As $m+1\leq k$ and $G$ is $k$-transitive, we have an element $\sigma\in G$ such that $\sigma(\beta_i)=\beta_i$ for all $1\leq i\leq m$ and $\sigma(\alpha)=\alpha'$. Consider $L=K(\beta_1,\beta_2, \dots,\beta_m)$ and $H=Gal(K_f/L)\subset G$. Clearly $\sigma \in H$. Also since $\sigma$ doesn't fix $\alpha$, we have $\alpha\not \in (K_f)^H=L$. But since $B$ is a minimal generating set, we have $L=K(B)=K_f$. Thus $\alpha\in L$ which gives a contradiction. \smallskip

 \item 

 Let $D=\{\delta_1,\delta_2,\dots, \delta_{m'}\} \subset R$ with $m'\leq k$. By a similar proof as in part (1) we have $\delta_i\not \in K(D\backslash \{\delta_i\})$ for all $1\leq i \leq m'$.\smallskip
 
 Now by Lemma \ref{minimal equivalent def} (3), there exists a minimal generating set $B'$ for $f$ over $K$. By part (1) we have $|B'|\geq k$. Consider any $m'$ many elements of $B'$ say $\beta_1,\beta_2, \dots, \beta_{m'}$. Since $G$ is $k$-transitive, we have a $\sigma\in G$ such that $\sigma(\beta_i)=\delta_i$ for all $1\leq i \leq m'$. Thus $D\subset B$ where $B=\sigma(B')=\{\sigma(\beta)\}_{\beta\in B'}$. We claim that $B$ is a minimal generating set for $f$ over $K$. Since $K(B')=K_f$ and $\sigma\in G$. Thus $K_f=\sigma(K_f)=\sigma(K(B'))=K(\sigma(B'))=K(B)$. Now let $A\subsetneq B$. Suppose $K(A)=K_f$. Let $\lambda=\sigma^{-1}$. Then as earlier $K_f=K(\lambda(A))$ where $\lambda(A)\subsetneq B'$. This is a contradiction as $B'$ is a minimal generating set. Thus $K(A)\neq K_f$ and we are done.

\smallskip

\item  
We have $k\geq 2$. Firstly we observe by part (1) and Proposition \ref{first prop} that $s_K(f)\geq k$. Now suppose on the contrary $r_K(f)\geq 2$. So we have two distinct roots of $f$ say $\alpha$ and $\alpha'$ such that $\alpha'\in K(\alpha)$. But by part (2) for $D=\{\alpha,\alpha'\}$ we get a contradiction. So $r_K(f)=1$. Also, we know from \cite{perlis2004roots} that $r_K(f)\cdot s_K(f)=n$. Thus indeed $s_K(f)=n$.

\smallskip

\item 
By part (1), the length of a shortest minimal cluster tower of $f$ is $\geq k$, thus we get by Proposition \ref{equivalence} (3) that the length of any cluster tower of $f$ is $\geq k$. \smallskip

Now we consider any cluster tower of $f$. \[ K \subseteq K(\beta_1) \subseteq K(\beta_1, \beta_2) \subseteq \dots \subseteq K(\beta_1, \beta_2, \ldots , \beta_s) = K_f.\]

We have to show that $[K(\beta_1,\beta_2,\dots, \beta_j):K]=\ ^nP_j$ for all $1\leq j\leq k$. It is enough to show that $[K(\beta_1,\beta_2,\dots, \beta_j):K(\beta_1,\beta_2,\dots, \beta_{j-1})]=n-j+1$ for all $1\leq j\leq k$. Thus we need to show that for each $1\leq j\leq k$, the polynomial $f_j(x)=f(x)/(\Pi_{i=1}^{j-1}(x-\beta_i))$ is irreducible over $L_{j-1}=K(\beta_1,\beta_2,\dots, \beta_{j-1})$.\smallskip

Suppose on the contrary that $f_j$ is reducible over $L_{j-1}$. Since $k<n$, we can consider two distinct roots of $f$ say $\alpha$ and $\alpha'$ which are respective roots of two distinct irreducible factors of $f_j$ over $L_{j-1}$ say $g$ and $g'$. As $j\leq k$ and $G$ is $k$-transitive, we have an element $\sigma\in G$ such that $\sigma(\beta_i)=\beta_i$ for all $1\leq i\leq j-1$ and $\sigma(\alpha)=\alpha'$. Consider $H=Gal(K_f/L_{j-1})\subset G$. Clearly $\sigma \in H$. Since $g$ is irreducible over $L_{j-1}$ and $g(\alpha)=0$, we have $\sigma(\alpha)=\alpha'$ to be a root of $g$ which gives a contradiction.\smallskip

In light of part (2), we can show that there are atleast $n(n-1)(n-2)\cdots (n-k+2)$ many distinct cluster towers by a similar argument as in proof of Proposition \ref{equivalence} (2).
\end{enumerate}

\smallskip

Now we assume parts (1)-(4) and establish that $G$ is $k$-transitive. We need to show that given any two ordered subsets of $R$ of size $k$, say $(\gamma_1,\gamma_2, \dots, \gamma_k)$ and $(\delta_1,\delta_2, \dots, \delta_k)$, we have $\sigma\in G$ such that $\sigma(\gamma_i)=\delta_i$ for all $1\leq i\leq k$. We will prove that for any $1\leq j\leq k$, there is a $\sigma_j\in G$ such that $\sigma_j(\gamma_i)=\delta_i$ for all $1\leq i\leq j$, by induction on $j\leq k$.\smallskip

The assertion is clear for $j=1$ as $G$ is a transitive subgroup of $\mathfrak{S}_n$. By induction assume that for $j<k$ the assertion is true. We will now prove for $j=k$. By part (2), $\{\delta_1,\delta_2,\dots, \delta_k\}$ is contained in a minimal generating set. We can consider the corresponding minimal cluster tower which is as follows. \[ K \subseteq K(\delta_1) \subseteq K(\delta_1, \delta_2) \subseteq \dots \subseteq K(\delta_1, \delta_2, \ldots , \delta_k) \subseteq \dots\]

By part (4) we have $[K(\delta_1, \delta_2,\dots, \delta_k):K(\delta_1,\delta_2,\dots, \delta_{k-1})]=n-k+1$. Thus the polynomial $f_k(x)=f(x)/(\Pi_{i=1}^{k-1}(x-\delta_i))$ is irreducible over $L_{j-1}=K(\delta_1,\delta_2,\dots, \delta_{k-1})$. Thus $H=Gal(K_f/L_{j-1})$ acts transitively on the set $R\backslash\{\delta_i\}_{i=1}^{k-1}$. By induction there is a $\sigma_{k-1}\in G$ such that $\sigma_{k-1}(\gamma_i)=\delta_i$ for all $1\leq i\leq k-1$. Let $\sigma_{k-1}(\gamma_k)=\delta$. Clearly $\delta\neq \delta_i$ for all $1\leq i\leq k-1$. As $\delta,\delta_k\in R\backslash\{\delta\}_{i=1}^{k-1}$, thus there is a $\lambda\in H$ such that $\lambda(\delta)=\delta_k$. Also as $\lambda\in H\subset G$, we have $\lambda(\delta_i)=\delta_i$ for all $1\leq i\leq k-1$. Now let $\sigma_j=\lambda\circ \sigma_{j-1}$. Thus $\sigma_j(\gamma_i)=\lambda(\sigma_{j-1}(\gamma_i))$. Hence for $1\leq i\leq k-1$, $\sigma_j(\gamma_i)=\lambda(\delta_i)=\delta_i$ and $\sigma_j(\gamma_k)=\lambda(\delta)=\delta_k$, so we are done.
\end{proof}

\begin{remark}
In the above, part (1) alone is not sufficient for $G$ to be $k$-transitive. Consider the case in \MainTheorem{} (1). Every minimal generating set for $f$ over $\mathbb{Q}$ has cardinality $\geq 2$. But $G={\rm Gal}(\mathbb{Q}_f/\mathbb{Q})=\mathbb{Z}/n \mathbb{Z} \rtimes (\mathbb{Z}/n \mathbb{Z})^{\times}$ is not $2$-transitive as there is no $\sigma\in G$ such that $\sigma(a)=a$ and $\sigma(a\zeta)=a\zeta^j$ where $(j,n)\neq 1$ (Note that $n$ is composite).
\end{remark}

\begin{remark}
    Theorem \ref{mult transtvty} (1), (2) and (4) are not true when $k=n$. In this case $G=\mathfrak{S_n}$ and by proof of \MainTheorem{} (2) for $k=n-1$, all minimal generating sets have cardinality $n-1$ (and minimal cluster towers have length $n-1$). 
\end{remark}

The following follows easily from Theorem \ref{mult transtvty} and its proof.

\begin{corollary}

    Consider an irreducible polynomial $f$ of degree $n>2$ over $K$ with the Galois group $G=Gal(K_f/K)$. Suppose $k+1<n$ and $G$ is $k$-transitive but not $(k+1)$-transitive. Then atleast one of the following is true.

\begin{enumerate}

\item There exists a subset $D$ of roots of $f$ of cardinality $k+1$ such that $D$ is not contained in any minimal generating set for $f$ over $K$, but for any $\delta\in D$, we have that $D\backslash \{\delta\}$ is contained in some minimal generating set for $f$ over $K$.\smallskip

\item There exists a cluster tower of $f$ of length $\geq k+1$ and the first $k$ terms in the degree sequence are $\ ^nP_1,\ ^nP_2 , \ldots,\ ^nP_{k}$\ but the $(k+1)$-th term is not $\ ^nP_{k+1}$\ .

    \end{enumerate}
\end{corollary}

By Theorem \ref{mult transtvty} (2) and transitivity of the Galois group we have,

\begin{corollary}
Consider an irreducible polynomial $f$ of degree $n$ over $K$. Every root of $f$ is contained in some minimal generating set for $f$ over $K$.
\end{corollary}



Now we give an application of Theorem \ref{mult transtvty}.



\begin{proof}[Proof of \MainTheorem{} (2) for $k=n-2$]

 By the results in \cite{brink2004alternating} on realizability of $\mathfrak{A}_n$ as a Galois group over hilbertian fields with characteristic $\neq 2$ and by results \cite{volklein1996groups} on hilbertian fields, we have ${\mathfrak A}_{n}$ to be realizable as a Galois group for infinitely many pairwise linearly disjoint Galois extensions over number field $K$. Thus by Proposition \ref{perm}, there exist infinitely many irreducible polynomials $f$ over $K$ of degree $n$ with Galois group ${\mathfrak A}_{n}$. Now we know that $\mathfrak{A}_n$ is $(n-2)$-transitive. Thus by Theorem \ref{mult transtvty} (1), any minimal generating set has cardinality $\geq n-2$. Now by Theorem \ref{mult transtvty} (4) and also noting that $|\mathfrak A_n|=n!/2=\ ^nP_{n-2}$ we have that length of any cluster tower is $n-2$. Thus every minimal generating set has cardinality $n-2$.\end{proof}

\begin{remark}

 Consider the case in proof above. 
 \begin{enumerate}
     \item The degree sequence of a cluster tower is independent of the ordering of the representatives of clusters of roots.\smallskip

     \item Any subset of roots of $f$ of cardinality $n-2$ is a minimal generating set and any minimal generating set is of this form. Hence there are $n\choose 2$ many minimal generating sets for $f$ over $K$.\smallskip

     \item  There are totally $n!/6$ many distinct cluster towers of $f$.
 \end{enumerate}
   
\end{remark}

Inspired from the discussion in \cite{gomez2011sharply}, we specialize Theorem \ref{mult transtvty} when the Galois group is sharply $k$-transitive.

\begin{theorem}\label{sharply mult trans}
     Consider an $n>2$ degree irreducible polynomial $f$ over $K$ with the Galois group $G=Gal(K_f/K)$. Suppose $k<n$. Then $G$ is sharply $k$-transitive if and only if we have all of the following.

     \smallskip
     \begin{enumerate}
    \item Every minimal generating set for $f$ over $K$ has cardinality $k$ and every subset of roots of $f$ of cardinality $k$ is a minimal generating set for $f$ over $K$. Thus there are exactly $n\choose k$ many distinct minimal generating sets.

 \smallskip

 \item The length of any cluster tower of $f$ is $k$ and the degree sequence is $(\ ^nP_1,\ ^nP_2 , \ldots,\ ^nP_{k})$\ . Thus there are exactly $\ ^nP_{k-1}$ many distinct cluster towers.
     \end{enumerate}
\end{theorem}

\begin{proof}
    We first assume that $G$ is sharply $k$-transitive and prove parts (1) and (2).\begin{enumerate}
\item  Let $D=\{\delta_1,\delta_2,\dots, \delta_k\}$ be a subset of roots of $f$ with $|D|=k$. Let $H=Gal(K_f/K(D))\subset G$. Suppose $\sigma\in H$, then $\sigma(\delta_i)=\delta_i$ for all $1\leq i\leq k$. Since $G$ is sharply $k$-transitive, $\sigma=id$. Thus $H=\{id\}$ and $K(D)=K_f$. By Theorem \ref{mult transtvty} (2), $D$ is a minimal generating set for $f$ over $K$.\smallskip

\item Follows from part (1) and Theorem \ref{mult transtvty} (4).

 \end{enumerate}

 Now we assume parts (1) and (2) and establish that $G$ is sharply $k$-transitive. Now parts (1) and (2) here imply the parts (2) and (4) in Theorem \ref{mult transtvty}. So by the proof of Theorem \ref{mult transtvty}, $G$ is $k$-transitive. If $G$ is not sharply $k$-transitive then there exists a subset $D$ of roots of $f$ of cardinality $k$ whose elements are fixed by a non-identity element. So we have $K(D)\neq K_f$ which contradicts part (1).\end{proof}

 \begin{remark}
   In the above case, the degree sequence of cluster tower is independent of the ordering of representatives of clusters of roots.
\end{remark}

 \begin{remark}
     Note that for $k+1<n$, $G$ being sharply $k$-transitive implies that $G$ is $k$-transitive but not $(k+1)$-transitive.
 \end{remark}

As a consequence of Theorem \ref{sharply mult trans}, we have the following.






\begin{proof}[Proof of \MainTheorem{} (3)] We use results (by Jordan, Dickson, Zassenhaus) on sharply $k$-transitive groups in \cite{dixon1996permutation} and the known cases of the famous Inverse Galois problem.

\begin{enumerate}[(\alph*)]

\item  The Mathieu groups $M_{12}$ and $M_{11}$ act transitively on 12 points and 11 points respectively and are sharply $5$-transitive and sharply $4$-transitive respectively. Now $M_{12}$ and $M_{11}$ are realizable as Galois groups over $\mathbb{Q}$ (Refer \cite{malle1999inverse}).

\smallskip

\item The group $PGL_2(\mathbb{F}_p)$ acts transitively on $p+1$ points and is sharply $3$-transitive. By the main result in \cite{arias2022locally}, $PGL_2(\mathbb{F}_p)$ is realizable as a Galois group over $\mathbb{Q}$.

\smallskip

\item The $1$-dimensional affine group $AGL(1,\mathbb{F}_q)\cong \mathbb{F}_q\rtimes \mathbb{F}_q^{\times}$ acts transitively on $q$ points and is sharply $2$-transitive. Clearly $AGL(1,\mathbb{F}_q)$ is solvable. By Shafarevich's theorem (\cite{Shafarevich1989}), $AGL(1,\mathbb{F}_q)$ is realizable as a Galois group over $\mathbb{Q}$. 
\end{enumerate}

Now applying Proposition \ref{perm}, we are done.\end{proof}

\begin{remark}
    Note that if for each group $PGL_2(\mathbb{F}_q)$ where $q>2$ is a prime power, there exists a finite Galois extension of the rationals with the Galois group equal to $PGL_2(\mathbb{F}_q)$, then we can generalize \MainTheorem{} (3)(b) for $n=q+1$ where $q>2$ is a prime power and $k=3$.
\end{remark}


\noindent {\it Acknowledgements:} The first author would like to acknowledge support of IIT Bombay Institute Post Doctoral Fellowship during the later part of this work and would like to thank his doctoral advisor Prof. Chandrasheel Bhagwat, IISER Pune for suggesting the problem of counting number of distinct cluster towers for a given polynomial under certain conditions. We also thank the referee for many pertinent comments and important suggestions.

\adresa
Shubham Jaiswal, Department of Mathematics, IIT Bombay, Powai, Mumbai 400 076, India,
e-mail: sjaiswal@math.iitb.ac.in

P Vanchinathan, e-mail: vanchinathan@gmail.com

\end{document}